# CONSISTENT ESTIMATES OF DEFORMED ISOTROPIC GAUSSIAN RANDOM FIELDS ON THE PLANE

By Ethan Anderes[1] and Sourav Chatterjee[2]

*University of California at Davis and University of California at Berkeley*

This paper proves fixed domain asymptotic results for estimating a smooth invertible transformation $f:\mathbb{R}^2 \to \mathbb{R}^2$ when observing the deformed random field $Z \circ f$ on a dense grid in a bounded, simply connected domain $\Omega$, where $Z$ is assumed to be an isotropic Gaussian random field on $\mathbb{R}^2$. The estimate $\hat{f}$ is constructed on a simply connected domain $U$, such that $\overline{U} \subset \Omega$ and is defined using kernel smoothed quadratic variations, Bergman projections and results from quasiconformal theory. We show, under mild assumptions on the random field $Z$ and the deformation $f$, that $\hat{f} \to R_\theta f + c$ uniformly on compact subsets of $U$ with probability one as the grid spacing goes to zero, where $R_\theta$ is an unidentifiable rotation and $c$ is an unidentifiable translation.

**1. Introduction.** The use of deformations to model nonstationary processes was first introduced to the spatial statistics literature by Sampson and Guttorp [37]. Their work, as well as that of subsequent authors (see, e.g., [13, 23, 34] and [38]) consider estimating the deformation $f$ when observing a deformed random field $Z \circ f$ at sparse observation locations with independent replicates of the random field.

Two recent papers, [10] and [5], study the different problem of estimating a deformation $f$ from dense observations of a single realization of a deformed isotropic random field $Z \circ f$ in two dimensions. These deformed isotropic random fields provide a flexible semi-parametric model of nonstationarity for random fields. In addition, this observation scenario is becoming increasingly important with the abundance of high resolution digital imagery and

Received June 2008; revised July 2008.
[1]Supported in part by NSF Postdoctoral Fellowship DMS-05-03227.
[2]Supported in part by a Sloan Research Fellowship in Mathematics and NSF Grant DMS-07-07054.

remote sensing. One of the more recent motivations for the deformation model under the one-realization observation scenario is gravitational lensing of the cosmic microwave background (CMB). The gravitational effect from intervening matter distort the CMB images to produce deformed random field observations. In the hope of improving estimates of cosmological parameters and the mass distribution in the universe (including dark matter) there is considerable interest in detecting and measuring the lensing of the CMB (see, e.g., [22] and [40]).

In this paper, we establish a strong consistency result for the estimation of the deformation $f$ when observing $Z \circ f$ on a dense grid in a bounded, simply connected domain in $\mathbb{R}^2$, as the grid spacing goes to zero. We first construct estimates of the complex dilatation and log-scale of the map $f$ (see Section 4), which converge uniformly on compact subsets of the observation region with probability one. Then, we construct a deformation estimate $\hat{f}$ on a subset of the observation region that converges uniformly on compact subsets with probability one. We show this result under mild assumptions on the map $f$ and the two dimensional isotropic random field $Z$.

Most attempts at recovering the deformation $f$ from a single realization of $Z \circ f$ rely on estimating local properties of $f$, usually related to the Jacobian of $f$, from the local behavior of the random field $Z \circ f$. Intuitively, the random field $Z \circ f$ is locally stretched and sheared by $f$, as determined by the Jacobian. One can clearly see the visual consequences of this shear, as seen in the left plot from Figure 1. When the random field $Z$ is isotropic, the identification of all four parameters of the Jacobian becomes difficult from the local behavior of $Z \circ f$. In particular, by decomposing the Jacobian matrix as $U\Lambda V^T$ (using singular value decomposition so that $U, V$ are orthogonal matrices, and $\Lambda$ is diagonal) the rotation matrix $U$ becomes particularly hard, if not impossible, to estimate when observing $Z \circ f$ in a small neighborhood. An important object for us is the complex dilatation and log-scale of $f$, determined by the Jacobian (and defined in Section 4), which are invariant under left multiplication of rotation matrices to the Jacobian. It is this invariance that allows us to estimate these parameters under the isotropy assumption for $Z$.

Guyon and Perrin [20] tackle the problem of developing consistent estimates of deformations in two dimensions and succeed in proving consistency within a subclass of deformations when observing random fields that are stationary but not isotropic. The subclass of deformations, however, is restrictive. In particular, it is not closed under post composition with rotations, which removes, in some sense, some complications for estimating the Jacobian generated by general deformations that can locally twist as well as stretch. On the other hand, Anderes and Stein [5] consider a large, nonparametric class of deformations. However, their results are methodological in nature and do not provide proofs of consistency. This paper contributes to bridging the gap between these two papers by considering the same flexible class of nonparametric deformations as in [5] while proving consistency for the estimated deformation as in [20]. For further references on densely observed deformed random fields see [9, 17, 29, 32, 33] and [35].

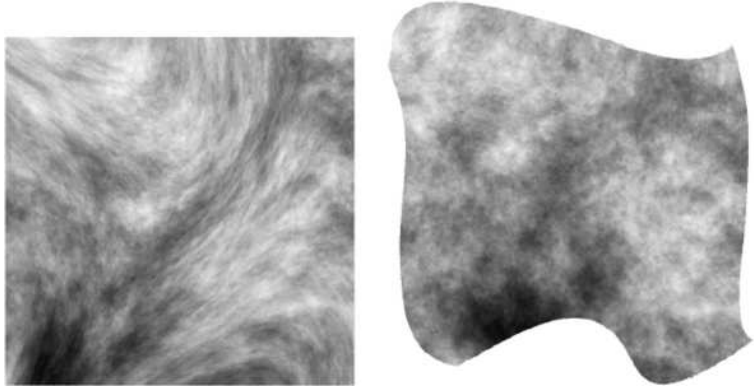

FIG. 1. *Left: One realization of a deformed isotropic random field. Right: The recovered isotropic random field using the estimated deformation from Figure 3.*

In this paper, we use kernel smoothed directional quadratic variations to estimate local properties of $f$, which are used to construct an estimate $\hat{f}$ of $f$. We establish sufficient conditions on the rate of bandwidth decay, in relation to the grid spacing, for strong uniform convergence of the kernel smoothed quadratic variations. There is a significant amount of literature studying the convergence of quadratic variations (see, e.g., [1, 6, 7, 8, 11, 15, 18, 19, 24, 25, 30] and [41]). One of the crucial inequalities used in many of the recent convergence results is the Hanson and Wright bound for quadratic forms [21]. Indeed, we also depend heavily on this bound, and we use it to establish Claim 2 in Section 3.2, which, in turn, gives uniform convergence on compact subsets for the estimated complex dilatation and log-scale and, ultimately, the convergence of the estimated deformation.

The kernel smoothed quadratic variations used in this paper are based on second order increments of the deformed process. Second order increments, rather than first order, are used in equation (4) to obtain sufficient spatial decorrelation for uniform convergence. Using higher order increments for quadratic variations is not new. They have been used in [24] and [7] for identification of a local fractional index and in [12] to identify the singularity function of a fractional process. The heuristic is that by increasing the order of increment, one can increase the rate of decay of the variance of the quadratic variation. However, this rate improvement holds only to a point, after which additional increments no longer improve the situation.

The interaction between the number of increments, the fractional index of the random field and the dimension of the domain of the random field is investigated in Chapter 3 of [3].

One of the main theoretical tools we use in this paper is the theory of quasiconformal maps. We believe this paper demonstrates how quasiconformal theory can provide a flexible theoretical framework for estimating smooth invertible transformations, whereby making these objects available to statisticians for modelling a diverse range of physical phenomena. In two dimensions, an important object in the theory of quasiconformal maps is the complex dilatation $\mu : \Omega \to \mathbb{D}$ (here $\Omega$ is the observation region for $Z \circ f$, and $\mathbb{D}$ is the unit disk in the complex plane). A more detailed discussion of the complex dilatation is presented in Section 4. Besides characterizing quasiconformal maps up to post composition with conformal maps, the complex dilatation $\mu$ has two other useful properties. First, $\mu$ can be interpreted as measuring the ellipticity and inclination of the local ellipse, which gets mapped to a local circle under the quasiconformal map that it characterizes. This is important to us for developing estimates of $\mu$ locally from $Z \circ f$. Second, the only requirement on $\mu$ is measurability and $\|\mu\|_\infty < 1$. In other words, it suffices to measurably assign eccentricity and inclinations of local infinitesimal ellipses, and, by keeping the eccentricity bounded, there is a quasiconformal map that sends these infinitesimal ellipses to circles (unique up to post composition with conformal maps). This property allows us to find a smooth invertible transformation that corresponds to the estimated complex dilatation.

The other object we use for estimating $f$ is the log-scale $\tau := \log |\partial f|$. We will discuss both $\mu$ and $\tau$ in detail later. However, it is worth while to notice that $\mu$ and $\tau$ provide enough information to uniquely specify the map $f$ up to an overall rotation and translation. As we will see, one of the difficulties with $\tau$, as compared to $\mu$, is that $\tau$ lies in a complicated subspace of functions mapping $\mathbb{R}^2$ to $\mathbb{R}$. This is where we employ the Bergman projection as a tool to overcome the restrictive nature of the log-scale parameter $\tau$.

Figure 2 illustrates the estimates $\hat{\mu}$ and $|\widehat{\partial f}|$ from the simulation shown in Figure 1. These are obtained by convolving the second order quadratic increments with a Gaussian smoothing kernel and transforming these smoothed increments as discussed in Section 4. An estimated deformation $\hat{f}$ corresponding to $\hat{\mu}$ and $|\widehat{\partial f}|$ is shown in Figure 3 (left) along with the true deformation (right). The image shown in Figure 1 (right) shows $Z \circ f \circ \hat{f}^{-1}$, which "unwinds" the deformed process in an attempt at recovering the isotropic process $Z$. Note that the deformation $\hat{f}$ in Figure 3 is constructed using methods from both this paper and from [5]. To be explicit, all of the estimation methods from this paper are used for $\hat{f}$ with the exception of the Bergman projection, where the methods from [5] were used. The reason for this is that the computational techniques are not yet developed, to the authors' knowledge, for accurate approximation of the Bergman projection used in constructing $\hat{f}$. However, since the Bergman space is a reproducing kernel Hilbert space, there is potential for accurate approximation using spline methodology.

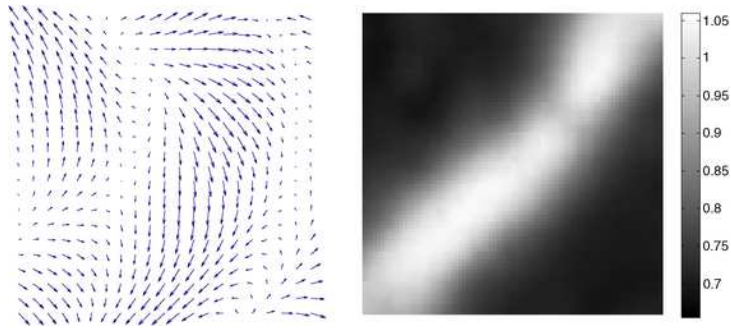

FIG. 2. *The estimated dilatation $\hat{\mu}$ (left) and scale $\widehat{|\partial f|}$ (right) using kernel smoothed second order directional increments and the results of Section 4.*

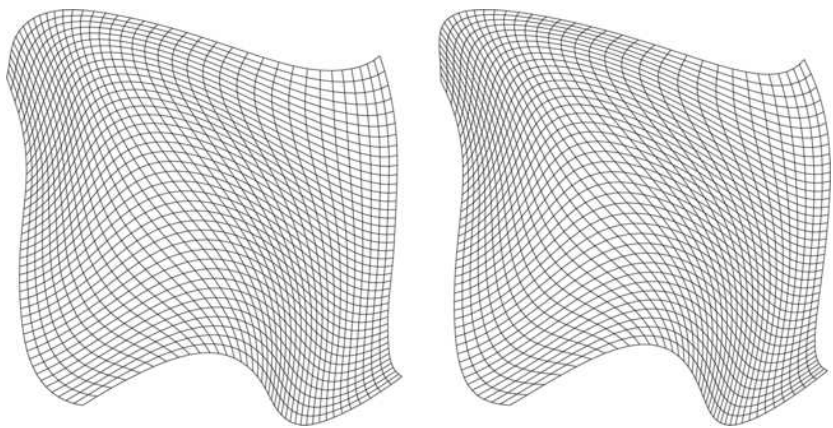

FIG. 3. *Left: The estimated deformation recovered from the estimated dilatation and scale shown in Figure 2. The deformation was constructed from the estimated dilatation and scale using methods outlined in [5]. Right: The true deformation $f$.*

There are three main parts to this paper. Section 2 discusses the assumptions on $Z$ and the smooth invertible map $f$. In Section 3, we show that kernel smoothed directional quadratic variations converge uniformly on compact subsets with probability one. These directional quadratic variations are then used, in Section 4, to get estimates $\hat{\mu}$ and $\hat{\tau} := \log \widehat{|\partial f|}$. Finally, in Section 5, we show how to convert $\hat{\mu}$ and $\hat{\tau}$ to an estimated map $\hat{f}$ on simply

connected subsets $U$ such that $\overline{U} \subset \Omega$, and we show that $\hat{f}$ converges to $f$ uniformly on compact subsets of $U$ with probability one.

**2. The random field $Z$ and the map $f$.** In this section, we list our assumptions on the isotropic random field $Z$ and the smooth invertible transformation $f$. This section starts with a brief discussion of our assumptions on the autocovariance function of $Z$. Then, a detailed discussion follows on our assumptions for the smooth invertible transformation $f : \mathbb{R}^2 \to \mathbb{R}^2$.

We require the following three conditions on the isotropic random field $Z$:

R1. $Z$ is a constant mean Gaussian process on $\mathbb{R}^2$ with autocovariance $R(|\mathbf{t} - \mathbf{s}|) = \text{cov}(Z(\mathbf{t}), Z(\mathbf{s}))$;

R2. $R(|t|) = R(0) - |t|^\alpha + o(|t|^{\alpha+\gamma})$, as $|t| \to 0$ for some $0 < \alpha < 2$, $\gamma > 0$;

R3. $R$ is $C^4$ away from the origin and there exists a $c > 0$ such that $|R^{(4)}(t)| \leq ct^{\alpha-4}$ for all sufficiently small $t > 0$.

The assumption R2 establishes the local quadratic variation behavior of the process $Z$ to be similar to that of a fractional Brownian sheet with Hurst index $\alpha/2$. Informally, the assumption R3 ensures that the second order increments of $Z$ have spatial decorrelation like that of a fractional Brownian sheet.

REMARK 1. Most of the following results can be extended, with some additional technical assumptions, to a larger class of autocovariance functions by replacing the principle term $|t|^\alpha$ in R2 with $L(|t|)|t|^\alpha$, where $L$ is a slowly varying function at 0. The main difference is that the quadratic variation process defined in equation (4) below will need to be normalized by $n^{-\alpha}L(1/n)$ instead of $n^{-\alpha}$.

REMARK 2. The class of autocovariances satisfying R2–R3 encompasses a broad range of random fields that are continuous but not differentiable. Examples include the Mátern autocovariance function with smoothness parameter $0 < \nu < 1$ (see [31, 39] and [42]) and the exponential family $\exp(-c|t|^\alpha)$, where $\alpha \in (0, 2)$. One way to extend our results to random fields with higher order differentiability is to use quadratic variations of higher order increments of the deformed random field to obtain sufficient spatial de-correlation. In the interest of space, we only prove results for the non differentiable case.

Our basic assumption on the smooth invertible map $f$ is that there exists a local affine approximation. In particular,

$$f(\mathbf{t} + \mathbf{h}) = f(\mathbf{t}) + J_f^\mathbf{t}\mathbf{h} + o(|\mathbf{h}|), \qquad (1)$$

where $J_f^{\mathbf{t}} := (\frac{\partial f_i}{\partial t_j}(\mathbf{t}))_{i,j}$ is the Jacobian of the map $f$ at $\mathbf{t}$. This local linear behavior is important, since we do not have replicates of the deformed random field $Z \circ f$, and, therefore, most of the statistical information is contained in the local variation of the process $Z \circ f$. When $f$ behaves locally like the Jacobian matrix transformation, the distribution of the random field $Z \circ f(\mathbf{t} + \mathbf{h})$, as $\mathbf{h}$ varies in a small neighborhood of the origin for a fixed $\mathbf{t}$, behaves similar to that of $Z(J_f^{\mathbf{t}} \mathbf{h})$. Therefore, one can hope to estimate parameters of the Jacobian $J_f^{\mathbf{t}}$ using the local quadratic variation of one realization of the process $Z \circ f$ near $\mathbf{t}$. Of course, higher order terms in a Taylor expansion may also be estimated; however, these presumably require smaller neighborhoods for accurate estimation.

In addition to the local affine behavior of $f$, we will need extra smoothness conditions. We discuss the following three notions of differentiability for a map $f: U \to V$ between planer open subsets $U, V$: Fréchet, Gâteaux and $C^1(U)$. For a point $\mathbf{x}_0 \in U$, $f$ is said to be Gâteaux differentiable at $\mathbf{x}_0$ in the direction $\mathbf{h}$ if the limit

$$\lim_{\varepsilon \to 0} \frac{f(\mathbf{x}_0 + \varepsilon \mathbf{h}) - f(\mathbf{x}_0)}{\varepsilon}$$

exists. A stronger notion of differentiability is Fréchet differentiable. The map $f$ is said to be Fréchet differentiable at $\mathbf{x}_0$ if there exits a continuous linear map $T: \mathbb{R}^2 \to \mathbb{R}^2$ such that

$$f(\mathbf{x}_0 + \mathbf{h}) - f(\mathbf{x}_0) = T(\mathbf{h}) + o(|\mathbf{h}|).$$

If such a $T$ exists, and $f$ has continuous partial derivatives, then $T$ is the map corresponding to left multiplication by the Jacobian matrix $J_f^{\mathbf{x}_0}$. Clearly, if a function $f$ is Fréchet differentiable at $\mathbf{x}_0$, then it is Gâteaux differentiable at $\mathbf{x}_0$ and the Gâteaux derivative in the direction $\mathbf{h}$ is equal to $J_f^{\mathbf{x}_0} \mathbf{h}$.

The third notion of differentiability $C^1(U)$ is satisfied if the partials $\frac{\partial f_i}{\partial x_j}$ exist and are continuous on $U$. This notion of differentiability is different, in that it is not defined pointwise. The reason for this is that there is not much one can say about the local behavior near $\mathbf{x}_0$ of a function where the partials exist. In fact, it may be neither Fréchet nor Gâteaux differentiable. If, however, we require the partials $\frac{\partial f_i}{\partial x_j}$ to exits and be continuous on $U$, then this is enough to imply Fréchet differentiability at all points in $U$. Even more, it is true that $f$ is $C^1(U)$, if and only if $\mathbf{x} \mapsto J_f^{\mathbf{x}}$ is continuous as a mapping from $U$ into the space of continuous linear functions on $\mathbb{R}^2$ (see Claim 8.9.1 from [14]).

Define the class of $C^1(U)$ diffeomorphisms to be the set of all continuous invertible maps $f: U \to \mathbb{R}^2$, such that $f$ is $C^1(U)$ and $f^{-1}$ is $C^1(V)$, where $V = f(U)$ is the range of $f$. By the Inverse Function Theorem, necessary and sufficient conditions are that $f$ be invertible, $C^1(U)$ and $\det J_f \neq 0$.

We write $C^1$ as short for $C^1(\mathbb{R}^2)$. Now, every $C^1$ diffeomorphism is Fréchet differentiable, and so there is a Jacobian matrix $J_f^{\mathbf{t}}$ such that

$$f(\mathbf{t}+\mathbf{h}) = f(\mathbf{t}) + J_f^{\mathbf{t}}\mathbf{h} + o(|\mathbf{h}|).$$

Moreover, the directional derivative in the direction $\theta$, denoted $\partial_\theta f$, is $J_f^{\mathbf{t}}\mathbf{u}_\theta$, where $\mathbf{u}_\theta = (\cos\theta, \sin\theta)$.

In the following paper, we will restrict the definition of $C^1$ diffeomorphisms to have $\det J_f > 0$ on $\mathbb{R}^2$, which characterizes the diffeomorphism to be orientation preserving. In some sense, this is a trivial restriction, since, when $f$ is a $C^1$ diffeomorphism, either $\det J_f > 0$ everywhere or $\det J_f < 0$ everywhere. These are referred to as orientation-preserving and orientation-reversing, respectively (see page 10 of [28]). Finally, define $C^r(U)$ diffeomorphisms to be the $C^1(U)$ diffeomorphisms with order $r$ continuous mixed partials.

We now list some consequences of a $C^2$ diffeomorphic assumption on $f$, which we use in the following proofs:

D1. $f$ is a quasiconformal map on bounded simply connected domains (see the Appendix of [4] for a definition of quasiconformal maps);
D2. $\sup_{\mathbf{t}\in\Theta} |\frac{f(\mathbf{t}+\varepsilon\mathbf{h})-f(\mathbf{t})}{\varepsilon} - J_f^{\mathbf{t}}\mathbf{h}| \to 0$ as $\varepsilon \to 0$ for every compact set $\Theta$;
D3. For any vector $\mathbf{h} \neq 0$ and compact set $\Theta$, there exists a constant $c$ such that $|\partial_{\mathbf{h}}^{(2,2)} R(|f(\mathbf{s}) - f(\mathbf{t})|)| \leq c|\mathbf{s}-\mathbf{t}|^{\alpha-4}$ for all $\mathbf{s},\mathbf{t}\in\Theta$ such that $\mathbf{s}\neq\mathbf{t}$ (note that $R$ and $\alpha$ are defined in assumptions R1–R3, and $\partial_{\mathbf{h}}^{(2,2)}$ is defined in the next section);
D4. For every compact subset $\Theta$, there exists constants $c_1, c_2 > 0$ such that $c_1|\mathbf{h}| \leq |J_f^{\mathbf{x}}\mathbf{h}| \leq c_2|\mathbf{h}|$ for all $\mathbf{h}$ and all $\mathbf{x}\in\Theta$;
D5. $|J_f^{\mathbf{t}}\mathbf{h}|^\alpha$ is Hölder continuous in $\mathbf{t}\in\Theta$ for any $\mathbf{h}$ and compact set $\Theta$.

Note that D3, D4 and D5 are the only statements that depend on the extra $C^2$ assumption rather than the $C^1$. The proofs of D1–D4 for $C^2$ diffeomorphisms are included in the Appendix of [4]. Notice that D5 follows from D4 and the fact that $J_f^{\mathbf{t}}$ has $C^1$ components.

**3. Kernel smoothed squared increments.** In this section, we study the convergence of the kernel smoothed squared increments of the deformed process $Y(\mathbf{x}) := Z \circ f(\mathbf{x})$ observed on some dense grid in a bounded, simply connected open subset $\Omega \subset \mathbb{R}^2$. The asymptotic regime we consider is as the grid spacing goes to zero and the region $\Omega$ stays fixed, which is sometimes called infill asymptotics.

For a fixed nonzero vector $\mathbf{h} \in \mathbb{R}^2$, let

$$\Delta_{\mathbf{h}} Y(\mathbf{t}) := Y(\mathbf{t}+\mathbf{h}) - Y(\mathbf{t}),$$
$$\Delta_{\mathbf{h}}^m Y(\mathbf{t}) := \Delta_{\mathbf{h}}\Delta_{\mathbf{h}}^{m-1} Y(\mathbf{t}).$$

If $\mathbf{t}$ is near $\partial\Omega$, computing $\Delta_\mathbf{h}^m$ may require observing $Y$ outside of $\Omega$. Therefore, we will suppose that we observe $Y$ on $\Omega$ plus some points within a small distance from the boundary $\partial\Omega$. Now, for a function of two variables $F(\mathbf{s},\mathbf{t})$, let $\Delta_\mathbf{h}^{(m,n)}F(\mathbf{s},\mathbf{t}):=\Delta_\mathbf{h}^m\Delta_\mathbf{h}^n F(\mathbf{s},\mathbf{t})$, where $\Delta_\mathbf{h}^m$ acts on the variable $\mathbf{s}$ and $\Delta_\mathbf{h}^n$ acts on the variable $\mathbf{t}$. Define $\partial_\mathbf{h}:=\mathbf{h}\cdot\nabla$ to be the directional derivative in the direction $\mathbf{h}$ and $\partial_\mathbf{h}^{(m,n)}F(\mathbf{s},\mathbf{t}):=\partial_\mathbf{h}^m\partial_\mathbf{h}^n F(\mathbf{s},\mathbf{t})$, where $\partial_\mathbf{h}^m$ acts on the variable $\mathbf{s}$ and $\partial_\mathbf{h}^n$ acts on $\mathbf{t}$. The following notation will be used throughout this paper:

$$(2) \qquad g(\mathbf{t}):=(8-2^{\alpha+1})|J_f^\mathbf{t}\mathbf{h}|^\alpha,$$

$$(3) \qquad \Omega_n:=\Omega\cap\{\mathbb{Z}^2/n\}.$$

That is, $\Omega_n$ is the grid of spacing $1/n$ in $\Omega$.

Here is a summary of the results of this section. First, we show Lemma 1, which establishes that $\mathsf{E}(\Delta_{\mathbf{h}/n}^2 Y(\mathbf{t}))^2\approx n^{-\alpha}g(\mathbf{t})$. Motivated by this lemma, we estimate $g(\mathbf{t})$ by locally averaging the squared increments $(\Delta_{\mathbf{h}/n}^2 Y(\mathbf{w}))^2/n^{-\alpha}$ for the points $\mathbf{w}$ near $\mathbf{t}$ then show that there is enough spatial decorrelation for convergence. To this end, define $B_{n,b}(\mathbf{t}):\Omega\to\mathbb{R}$ as

$$(4) \qquad B_{n,b}(\mathbf{t}):=\frac{1}{n^2b^2}\sum_{\mathbf{w}\in\Omega_n}K\left(\frac{\mathbf{w}-\mathbf{t}}{b}\right)\frac{(\Delta_{\mathbf{h}/n}^2 Y(\mathbf{w}))^2}{n^{-\alpha}}.$$

Here, $K$ is a convolution kernel satisfying certain conditions stated below. This section then culminates with Theorems 1 and 2 concerning the uniform convergence of $B_{n,b}$ and $\partial_\mathbf{u} B_{n,b}$ for some $\mathbf{u}\neq 0$. In particular, Theorem 1 shows that, under appropriate conditions,

$$\sup_{\mathbf{t}\in\Theta}|B_{n,b}(\mathbf{t})-g(\mathbf{t})|\longrightarrow 0 \qquad \text{w.p. } 1$$

as $n\to\infty$, $b\to 0$ and $n^{-1}b^{-3}\to 0$ for all compact sets $\Theta\subset\Omega$. Theorem 2 shows that, with some additional smoothness assumptions and $n^{-1}b^{-4}\to 0$, the directional derivatives $\partial_\mathbf{u} B_{n,b}$ converge uniformly on compact subsets w.p. 1 to $\partial_\mathbf{u} g$.

3.1. *Assumptions on $K$ and $\Omega$.* The assumptions on $K$ are as follows:

K1. $K$ has bounded, continuous first and second order mixed partial derivatives;
K2. $\int\int K(\mathbf{x})\,d\mathbf{x}=1$ and $\int\int|\mathbf{x}|K(\mathbf{x})\,d\mathbf{x}<\infty$.

Notice that these assumptions imply that $K$ is Riemann integrable, $K$ and the first partials of $K$ are Hölder continuous. Finally, we assume:

O1. $\Omega$ is a bounded simply connected domain of $\mathbb{R}^2$.

Notice that this assumption ensures that the number of points in $\Omega_n$ is of order $n^2$.

3.2. *Strong convergence of $B_{n,b}$ and $\partial_{\mathbf{u}} B_{n,b}$.* For the remainder of this section, let $Y := Z \circ f$, $B_{n,b}$ be defined as in (4), $g$ be defined as in (2), $\Omega_n$ be defined as in (3) and set $X_{\mathbf{t}} := B_{n,b}(\mathbf{t}) - \mathsf{E} B_{n,b}(\mathbf{t})$. In the following, 'universal constant' means any constant that does not depend on $n$, $b$, $\Theta$ or the process $(X_{\mathbf{t}})_{\mathbf{t} \in \Theta}$.

LEMMA 1. *Suppose R1–R2, O1 and $f$ is a $C^2$ diffeomorphism. Then,*
$$\sup_{\mathbf{t} \in \overline{\Omega}} \left| \frac{\mathsf{E}(\Delta^2_{\varepsilon \mathbf{h}} Y(\mathbf{t}))^2}{\varepsilon^\alpha} - g(\mathbf{t}) \right| \to 0$$
*as $\varepsilon \downarrow \infty$.*

PROOF. By assumption R2 we can write $R(|t|) = R(0) - |t|^\alpha + r(|t|)$, where $r(|t|) = o(|t|^{\alpha+\gamma})$ as $|t| \to 0$. Write $\mathsf{E}(\Delta^2_{\varepsilon \mathbf{h}} Y(\mathbf{t}_0))^2$ for $\mathbf{t}_0 \in \overline{\Omega}$ as a sum of two terms
$$\mathsf{E}(\Delta^2_{\varepsilon \mathbf{h}} Y(\mathbf{t}_0))^2 = \Delta^{(2,2)}_{\varepsilon \mathbf{h}} \{\mathrm{cov}(Y(\mathbf{s}), Y(\mathbf{t}))\}|_{\mathbf{s},\mathbf{t}=\mathbf{t}_0} = \mathcal{I}_1 + \mathcal{I}_2,$$
where
$$(5) \qquad \mathcal{I}_1 = \Delta^{(2,2)}_{\varepsilon \mathbf{h}} \{-|f(\mathbf{s}) - f(\mathbf{t})|^\alpha\}|_{\mathbf{s},\mathbf{t}=\mathbf{t}_0},$$
$$(6) \qquad \mathcal{I}_2 = \Delta^{(2,2)}_{\varepsilon \mathbf{h}} \{r(|f(\mathbf{s}) - f(\mathbf{t})|)\}|_{\mathbf{s},\mathbf{t}=\mathbf{t}_0}.$$
Write the increment operator $\Delta^2_{\varepsilon \mathbf{h}}$ as a linear filter, so that its action on a function $Q \colon \mathbb{R}^2 \to \mathbb{R}$ can be expressed as $\Delta^2_{\varepsilon \mathbf{h}} Q(\mathbf{t}) = \sum_{j=0}^{2} d_j Q(\mathbf{t} + \varepsilon \mathbf{s}_j)$, where $d_j = (-1)^{2-j} \binom{2}{j}$ and $\mathbf{s}_j = j\mathbf{h}$. Now, the first term can be computed as
$$(7) \qquad \mathcal{I}_1/\varepsilon^\alpha = -\sum_{i,j=0}^{2} d_i d_j |f(\mathbf{t}_0 + \varepsilon \mathbf{s}_i) - f(\mathbf{t}_0 + \varepsilon \mathbf{s}_j)|^\alpha / \varepsilon^\alpha$$
$$(8) \qquad = -\sum_{i,j=0}^{2} d_i d_j |J^{\mathbf{t}_0}_f (\mathbf{s}_i - \mathbf{s}_j)|^\alpha + o(1),$$
where $o(1) \to 0$ uniformly over $\mathbf{t}_0 \in \overline{\Omega}$ as $\varepsilon \to 0$ by D2 and D4.

Similarly, the second term can be computed as
$$(9) \qquad \mathcal{I}_2/\varepsilon^\alpha = \sum_{i,j=0}^{2} d_i d_j \frac{r(|f(\mathbf{t}_0 + \varepsilon \mathbf{s}_i) - f(\mathbf{t}_0 + \varepsilon \mathbf{s}_j)|)}{\varepsilon^\alpha},$$
where $\sup_{\mathbf{t}_0 \in \overline{\Omega}} |r(|f(\mathbf{t}_0 + \varepsilon \mathbf{s}_i) - f(\mathbf{t}_0 + \varepsilon \mathbf{s}_j)|)/\varepsilon^\alpha|$ converges to zero by R2, D2 and D4. Combining terms $\mathcal{I}_1$ and $\mathcal{I}_2$, we get
$$E(\Delta_{\varepsilon \mathbf{h}} Y(\mathbf{t}_0))^2/\varepsilon^\alpha \to -\sum_{i,j=0}^{2} d_i d_j |J^{\mathbf{t}_0}_f (\mathbf{s}_i - \mathbf{s}_j)|^\alpha$$

uniformly for $\mathbf{t}_0 \in \overline{\Omega}$. This completes the proof, since $-\sum_{i,j=0}^{2} d_i d_j |J_f^{\mathbf{t}_0}(\mathbf{s}_i - \mathbf{s}_j)|^\alpha = g(\mathbf{t}_0)$. □

LEMMA 2. *Suppose R1–R3, O1 and $f$ is a $C^2$ diffeomorphism. Then, there exists a constant $c > 0$ such that*

$$|\mathsf{E}\Delta^2_{\mathbf{h}/n} Y(\mathbf{t}) \Delta^2_{\mathbf{h}/n} Y(\mathbf{s})| \leq c n^{-4} |\mathbf{t} - \mathbf{s}|^{\alpha-4}$$

*for all $\mathbf{s}, \mathbf{t} \in \Omega$ such that $|\mathbf{s} - \mathbf{t}| > |3\mathbf{h}/n|$.*

PROOF. The idea is that

$$|\mathsf{E}\Delta^2_{\mathbf{h}/n} Y(\mathbf{t}) \Delta^2_{\mathbf{h}/n} Y(\mathbf{s})| = |\Delta^{(2,2)}_{\mathbf{h}/n} R(|f(\mathbf{s}) - f(\mathbf{t})|)|$$
$$\approx |n^{-4} \partial^{(2,2)}_{\mathbf{h}} R(|f(\mathbf{s}) - f(\mathbf{t})|)|$$
$$\leq c n^{-4} |\mathbf{s} - \mathbf{t}|^{\alpha-4}.$$

To make this precise, let $F(\mathbf{s}, \mathbf{t}) := \mathrm{cov}(Y(\mathbf{s}), Y(\mathbf{t})) = R(|f(\mathbf{s}) - f(\mathbf{t})|)$ and $H$ be the 2 by 2 matrix with each column $\mathbf{h}/n$. Then,

$$|\mathsf{E}\Delta^2_{\mathbf{h}/n} Y(\mathbf{t}) \Delta^2_{\mathbf{h}/n} Y(\mathbf{s})| = |\Delta^{(2,2)}_{\mathbf{h}/n} F(\mathbf{s}, \mathbf{t})|$$
$$= n^{-4} \left| \int_{[0,1]^2} \int_{[0,1]^2} \partial^{(2,2)}_{\mathbf{h}} F(\mathbf{s} + H\boldsymbol{\xi}, \mathbf{t} + H\boldsymbol{\eta}) \, d\boldsymbol{\xi} \, d\boldsymbol{\eta} \right|$$
$$\leq c_1 n^{-4} \int_{[0,1]^2} \int_{[0,1]^2} |\mathbf{s} - \mathbf{t} + H(\boldsymbol{\xi} - \boldsymbol{\eta})|^{\alpha-4} \, d\boldsymbol{\xi} \, d\boldsymbol{\eta}$$
$$\leq c_2 n^{-4} \sup_{\boldsymbol{\eta}, \boldsymbol{\xi} \in [0,1]^2} |\mathbf{s} - \mathbf{t} + \mathbf{h}(\xi_1 + \xi_2 - \eta_1 - \eta_2)/n|^{\alpha-4}$$
$$= c_2 n^{-4} \sup_{-1 \leq \tau \leq 1} |\mathbf{s} - \mathbf{t} + 2\mathbf{h}\tau/n|^{\alpha-4}$$
$$\leq c_3 n^{-4} |\mathbf{s} - \mathbf{t}|^{\alpha-4} \quad \text{when } |\mathbf{s} - \mathbf{t}| > |3\mathbf{h}/n|.$$

Notice that the above proof requires that $|\partial^{(2,2)}_{\mathbf{h}} F(\mathbf{s}, \mathbf{t})| \leq c_1 |\mathbf{s} - \mathbf{t}|^{\alpha-4}$ which is why we need D3. □

CLAIM 1. *Suppose R1–R3, K1–K2, O1 and $f$ is a $C^2$ diffeomorphism. Then,*

$$\sup_{\mathbf{t} \in \Theta} |\mathsf{E}(B_{n,b}(\mathbf{t})) - g(\mathbf{t})| \to 0$$

*for all compact subsets $\Theta \subset \Omega$ as $n \to \infty$, $b \to 0$ and $n^{-1} b^{-3} = o(1)$.*

PROOF. First, let $g_n(\mathbf{w}) := n^\alpha \mathsf{E}(\Delta^2_{\mathbf{h}/n} Y(\mathbf{w}))^2$ for $\mathbf{w} \in \Omega$. We show

$$\mathsf{E}(B_{n,b}(\mathbf{t})) = \frac{1}{n^2 b^2} \sum_{\mathbf{w} \in \Omega_n} K\left(\frac{\mathbf{w} - \mathbf{t}}{b}\right) g_n(\mathbf{w})$$

$$= \frac{1}{n^2 b^2} \sum_{\mathbf{w} \in \Omega_n} K\left(\frac{\mathbf{w} - \mathbf{t}}{b}\right) g(\mathbf{w}) + e_\mathrm{I}$$

$$= \frac{1}{b^2} \int \int_\Omega K\left(\frac{\mathbf{x} - \mathbf{t}}{b}\right) g(\mathbf{x}) \, d\mathbf{x} + e_\mathrm{II} + e_\mathrm{I}$$

$$= g(\mathbf{t}) + e_\mathrm{III} + e_\mathrm{II} + e_\mathrm{I},$$

where $e_\mathrm{I} = o(1)$, $e_\mathrm{II} = O(n^{-1} b^{-3})$ and $e_\mathrm{III} = o(1)$ uniformly on compact subsets of $\Omega$ as $n \to \infty$, $b \to 0$ and $n^{-1} b^{-3} \to 0$.

To show the results about $e_\mathrm{I}$ and $e_\mathrm{II}$, we need to control the error when approximating Riemann integrals of Hölder continuous functions on $\Omega$ by Riemann sums on $\Omega_n$. This error is bounded by the difference between the upper and lower Riemann sums, which is bounded by $a_\Omega a n^{-1}$, where $a$ is the Hölder constant of the function and $a_\Omega$ is a constant only depending on the region $\Omega$. To show $e_\mathrm{I} = o(1)$ and $e_\mathrm{II} = O(n^{-1} b^{-3})$ uniformly for $\mathbf{t} \in \Theta$, we will use the fact that $b^{-2} K((\cdot - \mathbf{t})/b)$ is Hölder continuous with Hölder constant $c b^{-3}$ for some constant $c$.

To show $e_\mathrm{I} = o(1)$, fix some compact subset $\Theta \subset \Omega$ and notice that

$$|e_\mathrm{I}| \leq \frac{1}{n^2 b^2} \sum_{\mathbf{w} \in \Omega_n} |K|\left(\frac{\mathbf{w} - \mathbf{t}}{b}\right) |g_n(\mathbf{w}) - g(\mathbf{w})|$$

$$\leq \left(\sup_{\mathbf{w} \in \overline{\Omega}} |g_n(\mathbf{w}) - g(\mathbf{w})|\right) \left(\frac{1}{n^2 b^2} \sum_{\mathbf{w} \in \Omega_n} |K|\left(\frac{\mathbf{w} - \mathbf{t}}{b}\right)\right).$$

The term $\sup_{\mathbf{w} \in \overline{\Omega}} |g_n(\mathbf{w}) - g(\mathbf{w})| \to 0$ by Lemma 1. Now, by the comments in the previous paragraph, the Riemann sum $\frac{1}{n^2 b^2} \sum_{\mathbf{w} \in \Omega_n} |K|(\frac{\mathbf{w}-\mathbf{t}}{b})$ is approximately $b^{-2} \int \int |K|((\mathbf{x} - \mathbf{t})/b) \, d\mathbf{x}$ (which is bounded) with error $O(n^{-1} b^{-3})$. Therefore, $e_\mathrm{I} \to 0$ as $n \to \infty$, $b \to 0$ and $n^{-1} b^{-3} \to 0$.

Similarly, to show $e_\mathrm{II} = O(n^{-1} b^{-3})$, we notice that the Hölder continuity of $K$ and $g$ are sufficient for the Riemann sums of $b^{-2} K((\cdot - \mathbf{t})/b) g(\cdot)$ to converge to the Riemann integral with an error $e_\mathrm{II} = O(n^{-1} b^{-3})$ uniformly in $\mathbf{t} \in \Theta$.

Finally, to show $e_\mathrm{III} = o(1)$, we need that

$$\frac{1}{b^2} \int \int_\Omega K\left(\frac{\mathbf{x} - \mathbf{t}}{b}\right) g(\mathbf{x}) \, d\mathbf{x} = \int \int_{(\Omega - \mathbf{t})/b} K(\mathbf{w}) g(b\mathbf{w} + \mathbf{t}) \, d\mathbf{w} \longrightarrow g(\mathbf{t})$$

as $b \to 0$ uniformly in $\mathbf{t} \in \Theta$. Here, we use the Hölder continuity of $g$ and assumption K2. Notice that the error term $e_\mathrm{III}$ does not converge to zero uniformly in $\mathbf{t} \in \overline{\Omega}$. This is why we can only show the result uniformly on compacts instead of uniformly on $\overline{\Omega}$. □

In what follows we will not only show convergence results about $B_{n,b}(\mathbf{t})$ but also $\partial_{\mathbf{u}} B_{n,b}(\mathbf{t})$, and $B_{n,b}(\mathbf{t}) - B_{n,b}(\mathbf{s})$ all of which have the form

$$(10) \quad Q_n = n^{\alpha-2} \sum_{\mathbf{i} \in \Omega_n} (\Delta^2_{\mathbf{h}/n} Y(\mathbf{i}))^2 \mathcal{K}_\mathbf{i},$$

where $\mathcal{K}_\mathbf{i}$ may depend on $\mathbf{t}, \mathbf{s} \in \Omega$ and the bandwidth parameter $b$.

CLAIM 2. *Let $Q_n$ be defined as in (10). Suppose R1–R3, O1 and $f$ is a $C^2$ diffeomorphism. In addition, suppose that there exists a function $G(\mathbf{t}, \mathbf{s}, b)$, such that $|\mathcal{K}_\mathbf{i}| \leq G(\mathbf{t}, \mathbf{s}, b)$ for all $\mathbf{s}, \mathbf{t} \in \Omega$ and $b$ in a neighborhood of the origin. Then, for all $\varepsilon > 0$, $b$ sufficiently small and $n$ sufficiently large,*

$$\mathsf{P}[|Q_n - \mathsf{E}Q_n| \geq \varepsilon] \leq c_1 \exp\left(-\frac{c_2 \varepsilon n^2}{G(\mathbf{t},\mathbf{s},b)} \wedge \frac{c_3 \varepsilon^2 n^2}{G(\mathbf{t},\mathbf{s},b)^2}\right),$$

*where $c_1, c_2, c_3$ are universal constants.*

PROOF. Writing $\mathcal{K}_\mathbf{i} = \mathcal{K}_\mathbf{i}^+ - \mathcal{K}_\mathbf{i}^-$, we get the decomposition $Q_n = Q_n^1 - Q_n^2$, where

$$Q_n^1 := n^{\alpha-2} \sum_{\mathbf{i} \in \Omega_n} (\Delta^2_{\mathbf{h}/n} Y(\mathbf{i}))^2 \mathcal{K}_\mathbf{i}^+,$$

$$Q_n^2 := n^{\alpha-2} \sum_{\mathbf{i} \in \Omega_n} (\Delta^2_{\mathbf{h}/n} Y(\mathbf{i}))^2 \mathcal{K}_\mathbf{i}^-.$$

Therefore,

$$\mathsf{P}(|Q_n - \mathsf{E}Q_n| \geq \varepsilon) \leq \mathsf{P}\left(|Q_n^1 - \mathsf{E}Q_n^1| \geq \frac{\varepsilon}{2}\right) + \mathsf{P}\left(|Q_n^2 - \mathsf{E}Q_n^2| \geq \frac{\varepsilon}{2}\right).$$

First, we find a bound for $\mathsf{P}(|Q_n^1 - \mathsf{E}Q_n^1| \geq \frac{\varepsilon}{2})$. Let $\Delta Y$ be the column vector with elements $n^{\alpha/2-1} \Delta^2_{\mathbf{h}/n} Y(\mathbf{i}) \mathcal{K}_\mathbf{i}^{+/2}$ for $\mathbf{i} \in \Omega_n$, where $\mathcal{K}_\mathbf{i}^{+/2} := (\mathcal{K}_\mathbf{i}^+)^{1/2}$. Let $\Sigma = \mathsf{E} \Delta Y \Delta Y^T$ be the covariance matrix for $\Delta Y$, so that $Q_n^1 = \Delta Y^T \Delta Y \stackrel{\mathcal{D}}{=} W^T \Sigma W$, where $W$ is a vector of i.i.d. standard Gaussian random variables. Also, let $\Sigma(\mathbf{i}, \mathbf{j})$ denote the matrix entries of $\Sigma$ for $\mathbf{i}, \mathbf{j} \in \Omega_n$.

Using the bound on quadratic forms for Gaussian random variables found in Hanson and Wright [21], we now get

$$\mathsf{P}[|W^T \Sigma W - \mathsf{E} W^T \Sigma W| \geq \varepsilon] \leq 2 \exp\left(-\frac{c_1 \varepsilon}{\|\Sigma_{\mathrm{abs}}\|_2} \wedge \frac{c_2 \varepsilon^2}{\|\Sigma_{\mathrm{abs}}\|_F^2}\right),$$

where $\|\cdot\|_2$ and $\|\cdot\|_F$ are the spectral and Frobenius matrix norms, respectively, and $\Sigma_{\text{abs}}$ is the matrix with elements $|\Sigma(\mathbf{i},\mathbf{j})|$ for $\mathbf{i},\mathbf{j}\in\Omega_n$. Now,

$$|\Sigma(\mathbf{i},\mathbf{j})| = \frac{\mathcal{K}_{\mathbf{i}}^{+/2}\mathcal{K}_{\mathbf{j}}^{+/2}}{n^{2-\alpha}}|\mathsf{E}\Delta_{\mathbf{h}/n}^2 Y(\mathbf{i})\Delta_{\mathbf{h}/n}^2 Y(\mathbf{j})|$$

$$\leq c_3 \frac{\mathcal{K}_{\mathbf{i}}^{+/2}\mathcal{K}_{\mathbf{j}}^{+/2}}{n^{2-\alpha}} n^{-4}|\mathbf{i}-\mathbf{j}|^{\alpha-4}$$

for all $|\mathbf{i}-\mathbf{j}|>|3\mathbf{h}/n|$ by Lemma 2. By assumption, $|\mathcal{K}_{\mathbf{i}}^{+/2}\mathcal{K}_{\mathbf{j}}^{+/2}|\leq G(\mathbf{t},\mathbf{s},b)$. Therefore, for all $\mathbf{i},\mathbf{j}\in\Omega_n$ such that $|\mathbf{i}-\mathbf{j}|>|3\mathbf{h}/n|$,

$$|\Sigma(\mathbf{i},\mathbf{j})| \leq c_3 G(\mathbf{t},\mathbf{s},b) n^{\alpha-6}|\mathbf{i}-\mathbf{j}|^{\alpha-4}. \tag{11}$$

To finish the proof, we show $\|\Sigma_{\text{abs}}\|_F^2 = O(|G(\mathbf{t},\mathbf{s},b)|^2 n^{-2})$ and $\|\Sigma_{\text{abs}}\|_2 = O(|G(\mathbf{t},\mathbf{s},b)|n^{-2})$ uniformly for all $\mathbf{t},\mathbf{s}\in\Omega$, $b$ sufficiently small and $n$ sufficiently large.

To show the bound for $\|\Sigma_{\text{abs}}\|_F^2$, notice that

$$\sum_{\mathbf{i},\mathbf{j}\in\Omega_n}\Sigma(\mathbf{i},\mathbf{j})^2 \leq \sum_{|\mathbf{i}-\mathbf{j}|\leq|3\mathbf{h}/n|}\Sigma(\mathbf{i},\mathbf{j})^2 + c_3^2 G^2 n^{2\alpha-8}\sum_{|\mathbf{i}-\mathbf{j}|>|3\mathbf{h}/n|} n^{-4}|\mathbf{i}-\mathbf{j}|^{2\alpha-8}$$

$$=^* \sum_{|\mathbf{i}-\mathbf{j}|\leq|3\mathbf{h}/n|}\Sigma(\mathbf{i},\mathbf{j})^2 + c_3^2 G^2 n^{2\alpha-8} O(n^{-2\alpha+6})$$

$$= O(n^{-2}G^2) + O(n^{-2}G^2),$$

where the last equality is because $|\Sigma(\mathbf{i},\mathbf{j})|\leq \max_{\mathbf{i}}\Sigma(\mathbf{i},\mathbf{i}) = O(n^{-2}G(\mathbf{t},\mathbf{s},b))$ by Lemma 1. To get $=^*$, there are some technical difficulties, but the heuristic is, when $0<\alpha<2$,

$$\sum_{|\mathbf{i}-\mathbf{j}|>|3\mathbf{h}/n|} n^{-4}|\mathbf{i}-\mathbf{j}|^{2\alpha-8} \asymp c_4 \int_{\{1/n<|\mathbf{x}|<1\}} |\mathbf{x}|^{2\alpha-8}\,d\mathbf{x}$$

$$= c_5 \int_{1/n}^1 r^{2\alpha-7}\,dr$$

$$= O(n^{-2\alpha+6}).$$

For the full details, see [3], Lemma 3, page 41. Finally, to show the bound for $\|\Sigma_{\text{abs}}\|_2$, notice that

$$\|\Sigma_{\text{abs}}\|_2 \leq \max_{\mathbf{i}\in\Omega_n}\sum_{\mathbf{j}\in\Omega_n}|\Sigma(\mathbf{i},\mathbf{j})|$$

$$\leq \max_{\mathbf{i}\in\Omega_n}\sum_{|\mathbf{i}-\mathbf{j}|\leq|3\mathbf{h}/n|}|\Sigma(\mathbf{i},\mathbf{j})| + c_3 G\max_{\mathbf{i}\in\Omega_n}\sum_{|\mathbf{i}-\mathbf{j}|>|3\mathbf{h}/n|}\frac{|\mathbf{i}-\mathbf{j}|^{\alpha-4}}{n^{6-\alpha}} \qquad \text{by (11)}$$

$$= O(n^{-2}G) + c_3 G n^{\alpha-4} \max_{\mathbf{i} \in \Omega_n} \sum_{|\mathbf{i}-\mathbf{j}|>|3\mathbf{h}/n|} n^{-2} |\mathbf{i} - \mathbf{j}|^{\alpha-4}$$

$$= O(n^{-2}G) + O(n^{-2}G),$$

where the last equality uses the fact that $0 < \alpha < 2$.

This establishes the desired bound for $\mathsf{P}(|Q_n^1 - \mathsf{E}Q_n^1| \geq \frac{\varepsilon}{2})$. The result for $\mathsf{P}(|Q_n^2 - \mathsf{E}Q_n^2| \geq \frac{\varepsilon}{2})$ is exactly similar. This completes the proof. □

COROLLARY 1. *Fix a point $\mathbf{t}_0 \in \Omega$, suppose R1–R3, O1 and let $f$ be a $C^2$ diffeomorphism. If $K$ is bounded, then*

$$|B_{n,b}(\mathbf{t}_0) - \mathsf{E}B_{n,b}(\mathbf{t}_0)| \xrightarrow{a.s.} 0$$

*as $n \to \infty$, $b \to 0$ and $n^{-1}b^{-2} = O(n^{-\beta})$ for some $\beta > 0$.*

PROOF. This follows by Claim 2 and Borel–Cantelli using $\mathcal{K}_{\mathbf{i}} = \frac{1}{b^2} K(\frac{\mathbf{i}-\mathbf{t}_0}{b})$ and $G(\mathbf{t}, \mathbf{s}, b) = b^{-2} \|K\|_\infty$. □

COROLLARY 2. *Fix a point $\mathbf{t}_0 \in \Omega$, suppose R1–R3, O1 and let $f$ be a $C^2$ diffeomorphism. If $K$ has continuous partial derivatives and $\mathbf{u} \neq 0$, then*

$$|\partial_{\mathbf{u}} B_{n,b}(\mathbf{t}_0) - \mathsf{E}\partial_{\mathbf{u}} B_{n,b}(\mathbf{t}_0)| \xrightarrow{a.s.} 0$$

*as $n \to \infty$, $b \to 0$ and $n^{-1}b^{-3} = O(n^{-\beta})$ for some $\beta > 0$.*

PROOF. This follows by Claim 2 and Borel–Cantelli using $\mathcal{K}_{\mathbf{i}} = \frac{1}{b^3} \times (\partial_{\mathbf{u}} K)(\frac{\mathbf{i}-\mathbf{t}_0}{b})$ and $G(\mathbf{t}, \mathbf{s}, b) = b^{-3} \|\partial_{\mathbf{u}} K\|_\infty$. □

The following corollary will be used for the uniform convergence of $B$ to $g$ in the next subsection. Remember that $X_{\mathbf{t}}$ is defined as $B_{n,b}(\mathbf{t}) - \mathsf{E}B_{n,b}(\mathbf{t})$.

COROLLARY 3. *Fix a point $\mathbf{t}_0 \in \Omega$, suppose R1–R3, O1 and let $f$ be a $C^2$ diffeomorphism. If $K$ is Hölder continuous, then*

$$(12) \qquad \mathsf{P}(|X_{\mathbf{t}} - X_{\mathbf{s}}| \geq \varepsilon) \leq c_1 \exp\left(-\frac{c_2 \varepsilon n^2 b^3}{|\mathbf{t}-\mathbf{s}|} \wedge \frac{c_3 \varepsilon^2 n^2 b^6}{|\mathbf{t}-\mathbf{s}|^2}\right),$$

*where $c_1, c_2, c_3$ are universal constants.*

PROOF. First, write $X_{\mathbf{t}} - X_{\mathbf{s}}$ in the form $Q_n - \mathsf{E}Q_n$, where $Q_n := B_{n,b}(\mathbf{t}) - B_{n,b}(\mathbf{s})$. Then, the corollary follows by Claim 2 using

$$\mathcal{K}_{\mathbf{i}} = \frac{1}{b^2} K\left(\frac{\mathbf{t}-\mathbf{i}}{b}\right) - \frac{1}{b^2} K\left(\frac{\mathbf{s}-\mathbf{i}}{b}\right),$$

so that $|\mathcal{K}_{\mathbf{i}}| \leq c b^{-3} |\mathbf{s} - \mathbf{t}|$ for some Hölder constant $c > 0$. □

3.3. *Uniform convergence of $B_{n,b}$ and $\partial_{\mathbf{u}} B_{n,b}$.* In this subsection, we use the results from the previous section to establish the uniform convergence of $B_{n,b}$ and the directional derivative $\partial_{\mathbf{u}} B_{n,b}$ on compact subsets of the observation region $\Omega$. These results are then used to establish consistent estimators of the complex dilatation $\mu$ and log-scale $\tau$ of the diffeomorphism $f$ in Section 4.

LEMMA 3. *For any $b > 0$ and $a \geq e$, we have*

$$\int_0^\infty (ae^{-bt^2} \wedge 1)\, dt \leq 2\sqrt{\frac{\log a}{b}} \quad \text{and} \quad \int_0^\infty (ae^{-bt} \wedge 1)\, dt \leq 2\frac{\log a}{b}.$$

PROOF. Given $b > 0$ and $a \geq e$, let $\gamma = \sqrt{b^{-1} \log a}$. Then,

$$\int_0^\infty (ae^{-bt^2} \wedge 1)\, dt = \gamma + \int_\gamma^\infty ae^{-bt^2}\, dt$$

$$\leq \gamma + \int_\gamma^\infty \frac{t}{\gamma} ae^{-bt^2}\, dt$$

$$= \gamma + \frac{1}{2b\gamma}$$

$$\leq 2\gamma \quad (\text{since } a \geq e \Rightarrow \gamma \geq 1/2b\gamma).$$

Similarly, putting $\nu = b^{-1} \log a$,

$$\int_0^\infty (ae^{-bt} \wedge 1)\, dt = \nu + \int_\nu^\infty ae^{-bt}\, dt = \nu + \frac{1}{b} \leq 2\nu.$$

This completes the proof of the lemma. $\square$

LEMMA 4. *Let $\Theta$ be a compact subset of $\Omega$. Suppose the assumptions in Corollary 3. Let $\mathbf{t}_1, \ldots, \mathbf{t}_m$ and $\mathbf{s}_1, \ldots, \mathbf{s}_m$ be arbitrary points in $\Theta$, where $m \geq 2$. Let $M = \max_{1 \leq i \leq m} |X_{\mathbf{t}_i} - X_{\mathbf{s}_i}|$ and $\delta = \max_{1 \leq i \leq m} |\mathbf{t}_i - \mathbf{s}_i|$. Then, for any $r > 0$,*

$$\mathsf{P}(M \geq r) \leq c_1 m \exp\left(-\frac{c_2 r n^2 b^3}{\delta}\right) + c_1 m \exp\left(-\frac{c_3 r^2 n^2 b^6}{\delta^2}\right),$$

*where $c_1$, $c_2$ and $c_3$ are universal constants.*

PROOF. By Corollary 3, for each $r > 0$ and $\mathbf{t}, \mathbf{s} \in \Theta$,

$$\mathsf{P}(|X_{\mathbf{t}} - X_{\mathbf{s}}| \geq r) \leq c_1 \exp\left(-\frac{c_2 r n^2 b^3}{|\mathbf{t} - \mathbf{s}|} \wedge \frac{c_3 r^2 n^2 b^6}{|\mathbf{t} - \mathbf{s}|^2}\right)$$

$$\leq c_1 \exp\left(-\frac{c_2 r n^2 b^3}{|\mathbf{t} - \mathbf{s}|}\right) + c_1 \exp\left(-\frac{c_3 r^2 n^2 b^6}{|\mathbf{t} - \mathbf{s}|^2}\right).$$

From the above bound, we see that

$$\mathsf{P}(M \geq r) \leq c_1 m \exp\left(-\frac{c_2 r n^2 b^3}{\delta}\right) + c_1 m \exp\left(-\frac{c_3 r^2 n^2 b^6}{\delta^2}\right).$$

This completes the proof. □

CLAIM 3. *Fix a compact set $\Theta \subset \Omega$ and a point $\mathbf{t}_0 \in \Theta$. Let $M = \max_{\mathbf{t} \in \Theta} |X_\mathbf{t} - X_{\mathbf{t}_0}|$ and suppose the assumptions in Corollary 3. Then, there exists universal constants $L_1$ and $L_2$ such that, for all $R > 0$, we have*

$$\mathsf{P}\left(M \geq \frac{R}{nb^3}\right) \leq L_1 \exp(-L_2 \min\{Rn, R^2\}).$$

PROOF. Suppose we have a sequence of finite sets $A_0, A_1, A_2, \ldots \subseteq \Theta$ and constants $c < 1 < B$ and $D$ satisfying the following properties:

(i) $A_0 = \{\mathbf{t}_0\}$;
(ii) $A_k \subseteq A_{k+1}$ for all $k$;
(iii) $|A_k| \leq B^k$ for all $k$, where $|A_k|$ denotes the cardinality of $A_k$;
(iv) For each $k \geq 1$ and each $\mathbf{t} \in A_k$ there exists a 'parent' $\mathbf{t}_p \in A_{k-1}$ such that $|\mathbf{t} - \mathbf{t}_p| \leq Dc^k$ (note that $\mathbf{t} = \mathbf{t}_p$ if $x \in A_{k-1}$);
(v) The sequence has a 'limiting denseness property' in the sense that, for any nonnegative continuous function $f$ on $\Theta$, we have $\max_{\mathbf{t} \in \Theta} f(\mathbf{t}) = \lim_{k \to \infty} \max_{\mathbf{t} \in A_k} f(\mathbf{t})$.

It is easy to see how the constants $c$, $B$ and $D$ can be chosen, with $D \propto \mathrm{diam}(\Theta)$, and the sets $\{A_k\}_{k \geq 1}$ constructed by successive dyadic partitioning. For each $k \geq 1$, let

$$M_k = \max_{\mathbf{t} \in A_k} |X_\mathbf{t} - X_{\mathbf{t}_p}|.$$

Applying the 'limiting denseness property,' we see that $M \leq \sum_{k=1}^\infty M_k$. For each $k \geq 1$, let $r_k = \frac{kc^k}{\sum_{j=1}^\infty jc^j}$. Then, $r_k > 0$ and $\sum_k r_k = 1$. Thus, for any $R > 0$,

$$\mathsf{P}\left(M \geq \frac{R}{nb^3}\right) \leq \mathsf{P}\left(M_k \geq \frac{Rr_k}{nb^3} \text{ for some } k \geq 1\right)$$
(13)
$$\leq \sum_{k=1}^\infty \mathsf{P}\left(M_k \geq \frac{Rr_k}{nb^3}\right).$$

By the tail bound in Lemma 4, we have

$$\mathsf{P}\left(M_k \geq \frac{Rr_k}{nb^3}\right) \leq c_1 B^k \exp\left(-\frac{c_2 Rr_k n}{Dc^k}\right) + c_1 B^k \exp\left(-\frac{c_3 R^2 r_k^2}{D^2 c^{2k}}\right)$$
$$\leq c_1 B^k (\exp(-c_5 Rkn) + \exp(-c_6 R^2 k^2)),$$

where $c_5$ and $c_6$ are universal constants. Thus, if $R > c_7 \log B$ for some suitably large constant $c_7$ (that need not depend on $n$), then, from the above bound and (13), we can conclude that

$$\mathsf{P}\bigg(M \geq \frac{R}{nb^3}\bigg) \leq c_8 \exp(-c_9 Rn) + c_8 \exp(-c_{10} R^2).$$

The condition $R > c_7 \log B$ can be dropped by choosing $c_8$ large enough, because $B$, $c_9$ and $c_{10}$ do not vary with $n$. This completes the proof. □

Now, we obtain Theorems 1 and 2, which are the main results of this section. These establish uniform convergence results for the kernel smoothed squared increments.

THEOREM 1. *Suppose* R1–R3, K1–K2, O1 *hold and $f$ is a $C^2$ diffeomorphism. Then, for all compact sets $\Theta \subset \Omega$,*

$$\sup_{\mathbf{t} \in \Theta} |B_{n,b}(\mathbf{t}) - g(\mathbf{t})| \longrightarrow 0 \qquad w.p.\ 1$$

*as $n \to \infty$, $b \to 0$ and $n^{-1}b^{-3} = O(n^{-\beta})$ for some $\beta > 0$, where $g(\mathbf{t}) := (8 - 2^{\alpha+1})|J_f^{\mathbf{t}}\mathbf{h}|^\alpha$.*

PROOF. Let $\Theta$ be a compact subset of $\Omega$ and let $\mathbf{t}_0 \in \Theta$. Let $b_n \to 0$ as $n \to \infty$, so that $n^{-1}b_n^{-3} \leq cn^{-\beta}$. Set $M_n = \sup_{\mathbf{t} \in \Theta} |X_\mathbf{t} - X_{\mathbf{t}_0}|$, where $X_\mathbf{t} := B_{n,b_n}(\mathbf{t}) - \mathsf{E} B_{n,b_n}(\mathbf{t})$. Then,

$$\sup_{\mathbf{t} \in \Theta} |B_{n,b_n}(\mathbf{t}) - g(\mathbf{t})| \leq \sup_{\mathbf{t} \in \Theta} |B_{n,b_n}(\mathbf{t}) - \mathsf{E} B_{n,b_n}(\mathbf{t})| + \sup_{\mathbf{t} \in \Theta} |\mathsf{E} B_{n,b_n}(\mathbf{t}) - g(\mathbf{t})|$$

$$\leq M_n + |B_{n,b_n}(\mathbf{t}_0) - \mathsf{E} B_{n,b_n}(\mathbf{t}_0)| + o(1) \qquad \text{by Claim 1.}$$

Now, by the bound in Claim 3, we have

$$\sum_{n=1}^{\infty} \mathsf{P}(M_n \geq n^{-\beta/2}) \leq \sum_{n=1}^{\infty} L_2 \exp(-L_3 \min\{n^{1+\beta/2}, n^\beta\}) < \infty.$$

Therefore, $M_n \to 0$ with probability one. Finally, $n^{-1}b_n^{-3} = O(n^{-\beta})$ for some $\beta > 0$ is a sufficient condition for $|B_{n,b_n}(\mathbf{t}_0) - \mathsf{E} B_{n,b_n}(\mathbf{t}_0)| \xrightarrow{\text{a.e.}} 0$ by Corollary 1. □

THEOREM 2. *Suppose the assumptions of Theorem 1 hold along with the additional assumption that $K$ is a compactly supported kernel and $f$ is a $C^3$ diffeomorphism. Then, for all compact sets $\Theta \subset \Omega$,*

$$\sup_{\mathbf{t} \in \Theta} |\partial_\mathbf{u} B_{n,b}(\mathbf{t}) - \partial_\mathbf{u} g(\mathbf{t})| \longrightarrow 0 \qquad w.p.\ 1$$

*as $n \to \infty$, $b \to 0$, $n^{-1}b^{-4} = O(n^{-\beta})$ for some $\beta > 0$ such that $\beta > 1 - 4\gamma$. Note that $\gamma > 0$ is the number appearing in assumption R2 and $g(\mathbf{t}) := (8 - 2^{\alpha+1})|J_f^{\mathbf{t}}\mathbf{h}|^\alpha$.*

PROOF. Let $\Theta$ be a compact subset of $\Omega$ and let $\mathbf{t}_0 \in \Theta$. Let $b_n \to 0$ as $n \to \infty$ so that $n^{-1} b_n^{-4} \leq c n^{-\beta}$. Set $M_n = \sup_{\mathbf{t} \in \Theta} |\tilde{X}_{\mathbf{t}} - \tilde{X}_{\mathbf{t}_0}|$, where $\tilde{X}_{\mathbf{t}} := \partial_{\mathbf{u}} B_{n,b_n}(\mathbf{t}) - \mathsf{E} \partial_{\mathbf{u}} B_{n,b_n}(\mathbf{t})$. Then,

$$\sup_{\mathbf{t} \in \Theta} |\partial_{\mathbf{u}} B_{n,b_n}(\mathbf{t}) - \partial_{\mathbf{u}} g(\mathbf{t})| \leq \sup_{\mathbf{t} \in \Theta} |\partial_{\mathbf{u}} B_{n,b_n}(\mathbf{t}) - \mathsf{E} \partial_{\mathbf{u}} B_{n,b_n}(\mathbf{t})|$$
$$+ \sup_{\mathbf{t} \in \Theta} |\mathsf{E} \partial_{\mathbf{u}} B_{n,b_n}(\mathbf{t}) - \partial_{\mathbf{u}} g(\mathbf{t})|$$
$$\leq M_n + |\partial_{\mathbf{u}} B_{n,b_n}(\mathbf{t}_0) - \mathsf{E} \partial_{\mathbf{u}} B_{n,b_n}(\mathbf{t}_0)| + o(1),$$

where $o(1)$ is established by similar methods to that of Claim 1 using the fact that $K$ is a $C^2$ compactly supported kernel.

Now, Corollary 2 establishes that $|\partial_{\mathbf{h}} B_{n,b_n}(\mathbf{t}_0) - \mathsf{E} \partial_{\mathbf{h}} B_{n,b_n}(\mathbf{t}_0)| \xrightarrow{\text{a.s.}} 0$. The proof that $M_n \xrightarrow{\text{a.s.}} 0$ is similar to Theorem 1, with the exception that, now, instead of the bound from Claim 3, one can show

$$\mathsf{P}\left(M_n \geq \frac{R}{nb^4}\right) \leq L_1 \exp(-L_2 \min\{Rn, R^2\}),$$

which is sufficient to prove the claim. □

**4. Consistent estimates of $(\mu, \tau)$.** At this point, it becomes advantageous to switch to complex variable notation, so that points in the plane $(x,y) \in \mathbb{R}^2$ correspond to points in the complex plane $x + iy \in \mathbb{C}$. Under this correspondence, $C^1$ diffeomorphisms $f(x,y) = (u(x,y), v(x,y))$ of $\mathbb{R}^2$ can be considered as $C^1$ diffeomorphisms of $\mathbb{C}$ by writing $f(x+iy) = u(x,y) + iv(x,y)$. For the remainder of the paper, we use $\operatorname{Re} f$ and $\operatorname{Im} f$ to denote the real and imaginary parts of the complex representation of the map $f$.

The main utility of switching to complex notation is that we can directly use the results and techniques of conformal and quasiconformal theory to establish consistent estimates of $f$. This section starts by defining the complex dilatation $\mu$ and log-scale $\tau$ of a $C^1$ diffeomorphism. We then conclude by showing how the results of the previous section can be used to construct consistent estimates of $\mu$ and $\tau$. In the next section, we show how the estimates of $\mu$ and $\tau$ can be used to establish consistent estimates of $f$.

For a function $f \in C^1(U)$, define the complex derivatives

$$\partial f := \frac{1}{2}\left(\frac{\partial f}{\partial x} - i \frac{\partial f}{\partial y}\right), \qquad \overline{\partial} f := \frac{1}{2}\left(\frac{\partial f}{\partial x} + i \frac{\partial f}{\partial y}\right).$$

The complex dilatation and the log-scale are then defined as

(14) $$\mu := \overline{\partial} f / \partial f,$$

(15) $$\tau := \log |\partial f|.$$

The complex dilatation $\mu$ characterizes the infinitesimal ellipse with inclination $\arg(-\mu/2)$ and eccentricity $\frac{1+|\mu|}{1-|\mu|}$ that gets mapped to an infinitesimal circle under the image of $f$. In addition, $\mu$ uniquely determines $f$ up to post composition with conformal maps. The log-scale $\tau$ is then used to recover the conformal post composition, so that, together, $\mu$ and $\tau$ uniquely determine $f$ up to a rotation and translation. For a short introduction to quasiconformal theory, see the appendices of [5] or [4]. For more a complete treatment, see [2, 26, 27] and [28].

In the previous section, we constructed a sequence of functions $B_{n,b_n}(\mathbf{t})$, which converge uniformly on compacts as $n \to \infty$ to $(8-2^{\alpha+1})|J_f^{\mathbf{t}}\mathbf{h}|^\alpha$, where $\mathbf{h} := (h_1, h_2)$ is a vector of our choice. Using complex variable notation, we can write $|J_f^{\mathbf{t}}\mathbf{h}| = |h\,\partial f + \overline{h}\,\overline{\partial} f|$, where $h = h_1 + ih_2$. By factoring out $\partial f$, we get $|h\,\partial f + \overline{h}\,\overline{\partial} f|^\alpha = |\partial f|^\alpha |h + \overline{h}\mu|^\alpha$. Now, by choosing $h = 1, i, 1+i$ (for increments in the north–south, east–west and diagonal directions), we can construct three functions $W_{1,n}$, $W_{2,n}$ and $W_{3,n}$, which converge to $|\partial f||1+\mu|$, $|\partial f||1-\mu|$ and $|\partial f||1+i+\mu(1-i)|$, respectively. In particular,

$$W_{1,n}(\mathbf{t}) \to \left|J_f^{\mathbf{t}}\begin{pmatrix}1\\0\end{pmatrix}\right|, \qquad W_{2,n}(\mathbf{t}) \to \left|J_f^{\mathbf{t}}\begin{pmatrix}0\\1\end{pmatrix}\right|, \qquad W_{3,n}(\mathbf{t}) \to \left|J_f^{\mathbf{t}}\begin{pmatrix}1\\1\end{pmatrix}\right|,$$

where the convergence is uniform in $\mathbf{t}$ on compacts subsets of $\Omega$ as $n \to \infty$. Notice that $W_{1,n}$, $W_{2,n}$ and $W_{3,n}$ are the factors of stretching that the affine transformation $J_f^{\mathbf{t}}$ applies to the lines in the horizontal, vertical and diagonal directions. Therefore, the points, $(W_{1,n}^{-1}, 0)$, $(0, W_{2,n}^{-1})$ and $(W_{3,n}^{-1}, W_{3,n}^{-1})$ asymptotically lie on the ellipse that gets mapped to a circle with unit radius under the affine transformation induced by the matrix $J_f^{\mathbf{t}}$. Since the general equation for an ellipse is $ax^2 + bxy + cy^2 = 1$, we have

$$a_n := W_{1,n}^2 \to a,$$
$$c_n := W_{2,n}^2 \to c,$$
$$b_n := W_{3,n}^2 - W_{1,n}^2 - W_{2,n}^2 \to b.$$

The area of the ellipse specified by $ax^2 + bxy + cy^2 = 1$ is $2\pi/\sqrt{4ac - b^2}$. Since $J_f^{\mathbf{t}}$ sends the ellipse $ax^2 + bxy + cy^2 = 1$ to the unit circle, we have that $\pi = \det(J_f^{\mathbf{t}})\frac{2\pi}{\sqrt{4ac-b^2}}$, which gives

$$\sqrt{4a_n c_n - b_n^2} \to 2\det(J_f^{\mathbf{t}}) = 2(|\partial f|^2 - |\overline{\partial} f|^2).$$

Now, to solve for $(\mu, \tau)$, first notice that $a_n + c_n \to |\partial f|^2(|1+\mu|^2 + |1-\mu|^2) = 2|\partial f|^2 + 2|\overline{\partial} f|^2$. Therefore, $\sqrt{4a_n c_n - b_n^2} + a_n + c_n \to 4|\partial f|^2$. Similarly,

(16) $\quad 4\operatorname{Re}(\mu) = |1+\mu|^2 - |1-\mu|^2,$

(17) $\quad 4\operatorname{Im}(\mu) = |1+i+\mu(1-i)|^2 - |1+\mu|^2 - |1-\mu|^2,$

which gives

$$\frac{a_n - c_n}{\sqrt{4a_nc_n - b_n^2} + a_n + c_n} \longrightarrow \mathrm{Re}(\mu), \tag{18}$$

$$\frac{b_n}{\sqrt{4a_nc_n - b_n^2} + a_n + c_n} \longrightarrow \mathrm{Im}(\mu), \tag{19}$$

$$\log(\sqrt{4a_nc_n - b_n^2} + a_n + c_n) \longrightarrow 2\tau + \log 4. \tag{20}$$

Therefore, under the conditions of Theorem 1, we can construct estimates $\hat{\mu}$ and $\hat{\tau}$ that converge to $\mu$ and $\tau$, respectively, where the convergence is uniform on compact subsets of $\Omega$ with probability one as $n \to \infty$. Moreover, under the extra conditions of Theorem 2, $\partial\hat{\mu} \to \partial\mu$ uniformly on compact subsets of $\Omega$ with probability one.

**5. Estimating $f$.** In this section, we show how to construct an estimate $\hat{f}$ on a simply connected domain $U$ such that $\overline{U} \subset \Omega$. Then, we show that $\hat{f}$ converges to $f$ uniformly on compact subsets of $U$. The construction of $\hat{f}$ is on $U$ instead of $\Omega$ because we need the uniform convergence of $\hat{\mu}$, $\partial\hat{\mu}$ and $\hat{\tau}$ to establish the convergence of $\hat{f}$. It is open as to whether one can construct $\hat{f}$ on the full observation region that converges uniformly on compact subsets.

We start, in Section 5.1, with a discussion on how $f$ can be recovered, uniquely up to a rotation and translation, from the true $\mu$ and $\tau$. This will indicate how we recover $\hat{f}$ from the estimated $\hat{\mu}$ and $\hat{\tau}$. Finally, we show $\hat{f} \to f$ uniformly on compact subsets of $U$ in Section 5.3.

5.1. *Recovering $f$ from $(\mu, \tau)$.* First, let $U$ be a simply connected domain such that $\overline{U} \subset \Omega$. The $C^1$ diffeomorphism $f$ now satisfies $f = g \circ f_\mu$ on $U$, where $f_\mu$ is the unique normalized quasiconformal map with dilatation $\mu$ that maps $U$ to the unit disk $\mathbb{D}$ (see the Appendix in [4]). Since $f$ and $f_\mu$ have the same complex dilatation, $g = f \circ f_\mu^{-1}$ is a conformal map defined on $\mathbb{D}$. It turns out that this decomposition of $f$ is useful, in that the complex dilatation $\mu$ determines $f_\mu$ and $\tau$ is used to recover the conformal map $g$.

To see how to recover $g$ from $\tau$, notice that $\partial f = \partial(g \circ f_\mu) = (g' \circ f_\mu)\partial f_\mu$. Therefore,

$$\log|g'| = \log|\partial f| \circ f_\mu^{-1} - \log|\partial f_\mu| \circ f_\mu^{-1}. \tag{21}$$

Since $g$ is conformal on $\mathbb{D}$, $\log g'$ is holomorphic on $\mathbb{D}$. Moreover, $\log g' = \log|g'| + i\arg(g')$. Therefore, using (21), $\mathrm{Re}\log g' = \tau \circ f_\mu^{-1} - \log|\partial f_\mu| \circ f_\mu^{-1}$, which can be recovered from $(\mu, \tau)$. Since the real and imaginary parts of holomorphic maps are harmonic conjugates, which are unique up to a constant, we can recover $\mathrm{Im}\log g' + \theta$, where $\theta \in \mathbb{R}$ is an unknown factor. Now, by exponentiating, we can recover $e^{i\theta}g'$. Then,

$$e^{i\theta}g(z) + c = \int_{z_0}^{z} e^{i\theta}g'(w)\,dw,$$

where the integral is taken over a line connecting $z_0$ to $z$. Therefore, $\mu$ and $\tau$ are sufficient to recover $f = g \circ f_\mu$ on $U$ up to a rotation and translation.

5.2. *Constructing $\hat{f}$ from $(\hat{\mu}, \hat{\tau})$.* The technique for recovering $\hat{f}$ in Section 5.1 would work if we knew that there existed a quasiconformal map $\hat{f}$ with complex dilatation $\hat{\mu}$ such that $\hat{\tau} = \log|\partial \hat{f}|$. Unfortunately, there is no simple condition on $\hat{\tau}$ for the existence of such an $\hat{f}$. The main problem is that we do not precisely measure $\log|g'|$, which is required to be harmonic. Instead, we only have an estimate of $\log|g'|$. The estimate, which is motivated by (21), is defined by

$$\widehat{\log|g'|} := \hat{\tau} \circ f_{\hat{\mu}}^{-1} - \log|\partial f_{\hat{\mu}}| \circ f_{\hat{\mu}}^{-1}. \tag{22}$$

Since $\widehat{\log|g'|}$ is not guaranteed to be harmonic, it may not always be possible to find the harmonic conjugate used to recover $\log g'$. In what follows, we notice that $\log g'$ is in the Bergman space of holomorphic functions with finite $L_2$ integrals. We then use the Bergman projection to find a holomorphic function whose real part approximates $\widehat{\log|g'|}$.

We define our estimate $\hat{f}$ of $f$ in the region $U$ as

$$\hat{f} = \hat{g} \circ f_{\hat{\mu}}, \tag{23}$$

where $f_{\hat{\mu}}$ is the unique normalized quasiconformal map sending $U$ to $\mathbb{D}$, and the function $\hat{g}$ is the holomorphic map, unique up to translations, defined on the unit disk $\mathbb{D}$ and satisfying

$$\hat{g}' = \exp(\mathcal{P}\widehat{\log|g'|}),$$

where the operator $\mathcal{P}$ is defined by

$$\mathcal{P}h(w) := \frac{2}{\pi} \int_\mathbb{D} \frac{h(z)}{(1-\bar{z}w)^2} \, dx\, dy - \operatorname{Re} h(0),$$

where $z = x + iy$. The integral transform in the above definition is the Bergman projection (see [16] for an introduction to Bergman spaces).

To motivate our choice of operator $\mathcal{P}$, we first mention that the true conformal map $g$ satisfies

$$(\mathcal{P}\log|g'|)(w) = \mathcal{P}\left(\frac{\log g' + \overline{\log g'}}{2}\right)(w)$$

$$= \frac{1}{\pi} \int_\mathbb{D} \frac{\log g'(z)}{(1-\bar{z}w)^2} \, dx\, dy + \frac{1}{\pi} \int_\mathbb{D} \frac{\overline{\log g'(z)}}{(1-\bar{z}w)^2} \, dx\, dy - \log|g'(0)|$$

$$= \log g'(w) + \overline{\log g'(0)} - \log|g'(0)|$$

$$= \log g'(w) + i\theta,$$

where $\theta = -\operatorname{Im}\log g'(0)$. In the above computation, we used the fact that for any conformal map $g$ defined on $\mathbb{D}$, the holomorphic function $\log g'$ is in the Bergman space $A^2(\mathbb{D})$ (this is true by Theorem 9.4 of [36] along with the fact that the Bloch space is a subset of the Bergman space). Therefore, the projection $\mathcal{P}$ can be used to recover the harmonic conjugate of $\log|g'|$ up to an unknown constant $\theta$. In what follows, we show $\widehat{\log|g'|} \to \log|g'|$ and $\mathcal{P}\widehat{\log|g'|} \to \mathcal{P}\log|g'|$ uniformly on compact subsets of $\mathbb{D}$ as $n \to \infty$.

### 5.3. $\hat{f}$ converges to $f$.

We show that, under appropriate conditions, $\hat{f}$ converges uniformly on compact subsets of $U$, with probability 1. First, we establish the following lemma.

LEMMA 5. *Suppose that $\mu_n, \mu \in C^2(U)$ are complex dilatations such that $\mu_n \overset{L_\infty(U)}{\to} \mu$ and $\partial\mu_n \overset{L_\infty(U)}{\to} \partial\mu$ on a bounded, simply connected domain $U$. Suppose, in addition, that one can extend $\mu_n, \mu$ to functions $\mu_n^*, \mu^* \in C^2(W)$ on a simply connected domain $W$ containing $\overline{U}$ such that $\mu_n \overset{L_\infty(W)}{\to} \mu$ and $\partial\mu_n \overset{L_\infty(W)}{\to} \partial\mu$. Then,*

$$\log|\partial f_{\mu_n}| \circ f_{\mu_n}^{-1} \to \log|\partial f_\mu| \circ f_\mu^{-1}, \tag{24}$$

$$\mathcal{P}\log|\partial f_{\mu_n}| \circ f_{\mu_n}^{-1} \to \mathcal{P}\log|\partial f_\mu| \circ f_\mu^{-1}, \tag{25}$$

*uniformly on compact subsets of $\mathbb{D}$ as $n \to \infty$.*

PROOF. First, decompose $f_{\mu_n}$ and $f_\mu$ so that,

$$f_{\mu_n} = h_n \circ \tilde{f}_n, \qquad f_\mu = h \circ \tilde{f},$$

where $\tilde{f}_n$ and $\tilde{f}$ are normalized quasiconformal maps on the whole plane with complex dilatations obtained by extending $\mu_n$ and $\mu$ to the whole plane by smoothly truncating to zero away from $U$. Here, $h_n$ and $h$ are conformal maps sending $\tilde{f}_n(U)$ and $\tilde{f}(U)$ to the unit disk $\mathbb{D}$, respectively. Note that the truncation must be done in such a way that $\tilde{\mu}_n$ has uniformly bounded compact support, $\tilde{\mu}_n$ is as smooth as $\mu_n$, $\tilde{\mu}_n \overset{L_\infty}{\to} \tilde{\mu}$ and $\partial\tilde{\mu}_n \overset{L_\infty}{\to} \partial\tilde{\mu}$ (existence guaranteed by the existence of the extensions $\mu^*, \mu_n^*$ on $W$). Now, the quasiconformal maps $\tilde{f}_n$ converge uniformly to $\tilde{f}$ on $U$, since $\tilde{\mu}_n \overset{L_\infty}{\to} \tilde{\mu}$ with uniformly bounded support (see Lemma 1 on page 55 of [2]).

Notice $\partial f_{\mu_n} = (h'_n \circ \tilde{f}_n)\partial\tilde{f}_n$, which gives

$$(\log\partial f_{\mu_n}) \circ f_{\mu_n}^{-1} = \log h'_n \circ h_n^{-1} + \log\partial\tilde{f}_n \circ f_{\mu_n}^{-1}. \tag{26}$$

See the Appendix in [4] for a discussion on how to define a continuous version of $\log\partial f_{\mu_n}$ and $\log\partial\tilde{f}_n$. To establish (24) and (25) we show that both terms in (26) converge uniformly on compact subsets of $\mathbb{D}$, as well as the result of applying the operator $\mathcal{P}$ to both terms.

For the first term, $\log h'_n \circ h_n^{-1}$, in (26), notice that $h_n^{-1} = \tilde{f}_n \circ f_{\mu_n}^{-1}$ and $h^{-1} = \tilde{f} \circ f_\mu^{-1}$. Since $\tilde{f}_n \to \tilde{f}$ uniformly on $U$, $\tilde{f}$ is Hölder continuous on $U$ and $f_{\mu_n}^{-1} \to f_\mu^{-1}$ uniformly on compact subsets of $\mathbb{D}$, $h_n^{-1} \to h^{-1}$ uniformly on compact subsets of $\mathbb{D}$ (this follows by Corollary 9 and Lemmas 10 and 13 of [4]). Since the functions $\xi_n := h_n^{-1}$ and $\xi := h^{-1}$ are conformal maps of $\mathbb{D}$, $\log \xi'_n$ and $\log \xi'$ are both holomorphic and $\log \xi'_n \to \log \xi'$ uniformly on compact subsets. Noticing that $\log \xi'_n = -\log h'_n \circ h_n^{-1}$ and $\log \xi' = -\log h' \circ h^{-1}$, gives

$$(27) \qquad \log h'_n \circ h_n^{-1} \to \log h' \circ h^{-1}$$

uniformly on compact subsets of $\mathbb{D}$ as $n \to \infty$. In addition, $\log \xi'_n$ and $\log \xi'$ are both in the Bergman space $A_2(\mathbb{D})$ (and are therefore unaffected by the Bergman projection), which establishes that

$$(28) \qquad \mathcal{P} \log h'_n \circ h_n^{-1} \to \mathcal{P} \log h' \circ h^{-1}$$

uniformly on compact subsets of $\mathbb{D}$ as $n \to \infty$.

For the second term $\log \partial \tilde{f}_n \circ f_{\mu_n}^{-1}$ in (26), notice that the results in the Appendix of [4] establish that $\log \partial \tilde{f}$ is Hölder continuous on $U$ and $\log \partial \tilde{f}_n \overset{L_\infty(U)}{\longrightarrow} \log \partial \tilde{f}$. Therefore, $\log \partial \tilde{f}_n \circ f_{\mu_n}^{-1}$ converges to $\log \partial \tilde{f} \circ f_\mu^{-1}$ uniformly on compact subsets. Moreover, since the continuity of $\log \partial \tilde{f}$ on $\mathbb{C}$ implies that it is bounded on $U$, $\log \partial \tilde{f}_n \circ f_{\mu_n}^{-1}$ also convergences in $L_2(\mathbb{D})$. Therefore,

$$(29) \qquad \log \partial \tilde{f}_n \circ f_{\mu_n}^{-1} \to \log \partial \tilde{f} \circ f_\mu^{-1},$$

$$(30) \qquad \mathcal{P} \log \partial \tilde{f}_n \circ f_{\mu_n}^{-1} \to \mathcal{P} \log \partial \tilde{f} \circ f_\mu^{-1}$$

uniformly on compact subsets of $\mathbb{D}$ as $n \to \infty$. The last convergence is due to the fact that the Bergman projection is a bounded operator on $L_2(\mathbb{D})$ and that convergence in the Bergman space implies convergence on compacts. Finally, (27), (28), (29) and (30) establishes the lemma. □

THEOREM 3. *Suppose* R1–R3*,* O1*,* K1–K2*,* K *is a compactly supported kernel,* $f$ *is a* $C^3$ *diffeomorphism,* $n \to \infty$*,* $b \to 0$ *and* $n^{-1}b^{-4} = O(n^{-\beta})$ *for some* $\beta > 0$ *such that* $\beta > 1 - 4\gamma$*. Let* $U$ *be a simply connected open subset of the observation region* $\Omega$ *such that* $\overline{U} \subset \Omega$*. Then, the estimated map* $\hat{f}$*, defined on* $U$ *by* (23)*, converges to* $e^{i\theta} f + c$ *uniformly on compact subsets of* $U$ *with probability one, where* $\theta$ *is an unidentifiable rotation angle and* $c$ *is an unidentifiable translation.*

PROOF. The results of Theorems 1 and 2, along with the comments made in Section 4, establish that $\hat{\mu} \to \mu$, $\partial \hat{\mu} \to \partial \mu$ and $\hat{\tau} \to \tau$ uniformly on $U$ with probability one as $n \to \infty$. Therefore, $f_{\hat{\mu}} \to f_\mu$ uniformly on compact subsets of $U$ (by Corollary 9 of [4]). It is now sufficient to show that $\hat{g}$ converges uniformly on compact subsets of $\mathbb{D}$ to $g$ (where sufficiency is by Lemma 11 of [4]).

We first show

$$\widehat{\log |g'|} \to \log |g'|, \tag{31}$$

$$\mathcal{P}\widehat{\log |g'|} \to \mathcal{P} \log |g'| \tag{32}$$

uniformly on compact subsets of $\mathbb{D}$ as $n \to \infty$. Remember, $\widehat{\log |g'|}$ is defined by

$$\widehat{\log |g'|} := \hat{\tau} \circ f_{\hat{\mu}}^{-1} - \log |\partial f_{\hat{\mu}}| \circ f_{\hat{\mu}}^{-1}. \tag{33}$$

Lemma 5 immediately establishes the required convergence for the second term $\log |\partial f_{\hat{\mu}}| \circ f_{\hat{\mu}}^{-1}$. The first term $\hat{\tau} \circ f_{\hat{\mu}}^{-1}$ converges to $\tau \circ f_\mu^{-1}$ both uniformly on compacts of $\mathbb{D}$ and in $L_2(\mathbb{D})$. This follows by Lemma 13 of [4], since $\hat{\tau} \xrightarrow{L_\infty(U)} \tau$ and that $\tau$ is Hölder continuous and bounded on $U$ (since $f$ is assumed to be a $C^3$ diffeomorphism). Since the Bergman projection is a bounded operator from $L_2(\mathbb{D})$ to $A_2(\mathbb{D})$, we also have that $\mathcal{P}\hat{\tau} \circ f_{\hat{\mu}}^{-1} \xrightarrow{A_2(\mathbb{D})} \mathcal{P}\tau \circ f_\mu^{-1}$. The convergence is also uniformly on compact subsets by Lemma 12 of [4]. This establishes (31) and (32).

Now, by the comments made in Section 5.2, $\mathcal{P} \log |g'| = \log g' + i\theta$. Therefore, by (31) and (32), $\exp(\mathcal{P}\widehat{\log |g'|})$ converges uniformly on compacts to $e^{i\theta} g'$. Since $U$ is simply connected, $\exp(\mathcal{P}\widehat{\log |g'|})$ has an antiderivative $\hat{g}$ such that $\hat{g}' = \exp(\mathcal{P}\widehat{\log |g'|})$ and $\hat{g} \to e^{i\theta} g + c$ uniformly on compact subsets of $\mathbb{D}$. Therefore, $\hat{f}$ converges to $e^{i\theta} f + c$ uniformly on compact subsets of $U$. □

**Acknowledgments.** The authors would like to thank the editor and three referees for many useful comments and suggestions.

## REFERENCES

[1] ADLER, R. J. and PYKE, R. (1993). Uniform quadratic variation for Gaussian processes. *Stochastic Process. Appl.* **48** 191–209. MR1244542
[2] AHLFORS, L. V. (2006). *Lectures on Quasiconformal Mappings. University Lecture Series* **38** Amer. Math. Soc., Providence, RI. MR2241787
[3] ANDERES, E. B. (2005). Estimating deformations of isotropic Gaussian random fields. Ph.D. thesis, Univ. Chicago.

[4] Anderes, E. B. and Chatterjee, S. (2008). Consistent estimates of deformed isotropic Gaussian random fields on the plane. Technical Report 739, Statistics Dept., Univ. California at Berkeley. Available at http://www.stat.berkeley.edu/tech-reports/739.pdf.

[5] Anderes, E. B. and Stein, M. L. (2008). Estimating deformations of isotropic Gaussian random fields on the plane. *Ann. Statist.* **36** 719–741. MR2396813

[6] Baxter, G. (1956). A strong limit theorem for Gaussian processes. *Proc. Amer. Math. Soc.* **7** 522–527. MR0090920

[7] Benassi, A., Cohen, S., Istas, J. and Jaffard, S. (1998). Identification of filtered white noises. *Stochastic Process. Appl.* **75** 31–49. MR1629014

[8] Berman, S. M. (1967). A version of the Lévy–Baxter theorem for the increments of Brownian motion of several parameters. *Proc. Amer. Math. Soc.* **18** 1051–1055. MR0222958

[9] Clerc, M. and Mallat, S. (2002). The texture gradient equation for recovering shape from texture. *IEEE Trans. on Pattern Analysis and Machine Intelligence* **24** 536–549.

[10] Clerc, M. and Mallat, S. (2003). Estimating deformations of stationary processes. *Ann. Statist.* **31** 1772–1821. MR2036390

[11] Cohen, S., Guyon, X., Perrin, O. and Pontier, M. (2006). Identification of an isometric transformation of the standard Brownian sheet. *J. Statist. Plann. Inference* **136** 1317–1330. MR2253765

[12] Cohen, S., Guyon, X., Perrin, O. and Pontier, M. (2006). Singularity functions for fractional processes: Application to the fractional Brownian sheet. *Ann. Inst. H. Poincaré. Probab. Statist.* **42** 187–205. MR2199797

[13] Damian, D., Sampson, P. and Guttorp, P. (2001). Bayesian estimation of semiparametric non-stationary spatial covariance structures. *Environmetrics* **12** 161–178.

[14] Dieudonné, J. (1960). *Foundations of Modern Analysis*. Academic Press, New York. MR0120319

[15] Dudley, R. M. (1973). Sample functions of the Gaussian process. *Ann. Probab.* **1** 66–103. MR0346884

[16] Duren, P. and Schuster, A. (2004). *Bergman Spaces. Mathematical Surveys and Monographs* **100**. Amer. Math. Soc., Providence, RI. MR2033762

[17] Gårding, J. (1992). Shape from texture for smooth curved surfaces in perspective projection. *J. Math. Imaging Vision* **2** 327–350.

[18] Gladyshev, E. G. (1961). A new limit theorem for stochastic processes with Gaussian increments. *Theory Probab. Appl.* **6** 52–61. MR0145574

[19] Guyon, X. and Leon, G. (1989). Convergence en loi des h-variations d'un processus gaussien stationnaire. *Ann. Inst. H. Poincaré* **25** 265–282. MR1023952

[20] Guyon, X. and Perrin, O. (2000). Identification of space deformation using linear and superficial quadratic variations. *Statist. Probab. Lett.* **47** 307–316. MR1747492

[21] Hanson, D. L. and Wright, F. T. (1971). A bound on tail probabilities for quadratic form in independent random variables. *Ann. Math. Statist.* **42** 1079–1083. MR0279864

[22] Hu, W. (2001). Dark synergy: Gravitational lensing and the cmb. *Phys. Rev. D* **65**.

[23] Iovleff, S. and Perrin, O. (2004). Estimating a nonstationary spatial structure using simulated annealing. *J. Comput. Graph. Statist.* **13** 90–105. MR2044872

[24] Istas, J. and Lang, G. (1997). Quadratic variations and estimation of the local Hölder index of a Gaussian process. *Ann. Inst. H. Poincaré Probab. Statist.* **33** 407–436. MR1465796

[25] Klein, R. and Gine, E. (1975). On quadratic variation of processes with Gaussian increments. *Ann. Probab.* **3** 716–721. MR0378070

[26] Krushkal', S. L. (1979). *Quasiconformal Mappings and Riemann Surfaces*. V. H. Winston & Sons, Washington, DC. MR0536488

[27] Ławrynowicz, J. (1983). *Quasiconformal Mappings in the Plane: Parametrical Methods. Lecture Notes in Mathematics* **978**. Springer, Berlin. MR0702025

[28] Lehto, O. and Virtanen, K. I. (1965). *Quasiconformal Mappings in the Plane*. Springer, New York. MR0344463

[29] Leon, J. and Ortega, J. (1989). Weak convergence of different types of variation for biparametric Gaussian processes. In *Limit Theorems in Probability and Statistics. Colloqnia Mathematica Societatis János Bolyali* **57** 349–364.

[30] Lévy, P. (1940). Le mouvement brownien plan. *Amer. J. Math.* **62** 487–550. MR0002734

[31] Loh, W. (2005). Fixed-domain asymptotics for a subclass of matern-type Gaussian random fields. *Ann. Statist.* **33** 2344–2394. MR2211089

[32] Malik, J. and Rosenholtz, R. (1997). Computing local surface orientation and shape from texture for curved surfaces. *Int. J. Comp. Vision* **23** 149–168.

[33] Perrin, O. (1998). Functional convergence in distribution of quadratic variations for a large class of Gaussian processes: Application to a time deformation model. Technical report, Unit of Biometrics at Avignon.

[34] Perrin, O. and Meiring, W. (1999). Identifiability for non-stationary spatial structure. *J. Appl. Probab.* **36** 1244–1250. MR1746409

[35] Perrin, O. and Senoussi, R. (2000). Reducing non-stationary random fields to stationarity and isotropy using a space deformation. *Statist. Probab. Lett.* **48** 23–32. MR1767607

[36] Pommerenke, C. (1975). *Univalent Functions*. Vandenhoeck & Ruprecht, Göttingen, Germany. MR0507768

[37] Sampson, P. and Guttorp, P. (1992). Nonparametric estimation of nonstationary spatial covariance structure. *J. Amer. Statist. Assoc.* **87** 108–119.

[38] Schmidt, A. and O'Hagan, A. (2003). Bayesian inference for nonstationary spatial covariance structure via spatial deformations. *J. Roy. Statist. Soc. Ser. B* **65** 745–758.

[39] Stein, M. L. (1999). *Interpolation of Spatial Data: Some Theory for Kriging*. Springer, New York. MR1697409

[40] Stompor, R. and Efstathiou, G. (1999). Gravitational lensing of the cosmic microwave background anisotropies and cosmological parameter estimation. *Monthly Notices of the Royal Astronomical Society* **302** 735–747.

[41] Strait, P. T. (1969). On Berman's version of the Lévy–Baxter theorem. *Proc. Amer. Math. Soc.* **23** 91–93. MR0246358

[42] Zhang, H. (2004). Inconsistent estimation and asymptotically equal interpolations in model-based geostatistics. *J. Amer. Statist. Assoc.* **99** 250–261. MR2054303

DEPARTMENT OF STATISTICS
UNIVERSITY OF CALIFORNIA AT DAVIS
4214 MATHEMATICAL SCIENCES BLDG.
DAVIS, CALIFORNIA 95616
USA
E-MAIL: anderes@stat.ucdavis.edu

DEPARTMENT OF STATISTICS
UNIVERSITY OF CALIFORNIA AT BERKELEY
367 EVANS HALL
BERKELEY, CALIFORNIA 94720-3860
USA
E-MAIL: sourav@stat.berkeley.edu